\documentclass[11pt]{article}
\usepackage{amsmath,amssymb}
\usepackage[french]{babel}

\catcode`\@=12

\def\A{{\cal A}}
\def\B{{\cal B}}
\def\C{{\cal C}}
\def\Ker{\mathrm {Ker}}
\def\Im{\mathrm {Im}}
\let\la=\lambda

\let\til=\widetilde
\def\scal(#1,#2){\langle#1,#2\rangle}
\newtheorem{theoreme}{Théorème}

\title{Facteurs disjoints des transformations m\'elangeantes.}
\author{Fran\c{c}ois Parreau}
\date{}
\begin{document}
\maketitle
\begin{abstract}
We show that any non-mixing automorphism of a standard probability space has a factor disjoint from all mixing automorphism.
\end{abstract}
Soient $T$ et $S$ deux automorphismes de deux espaces de probabilit\'e
standard $(X,\B,\mu)$ et $(Y,\C,\nu)$. \'Etant donn\'e un couplage non
trivial $\la$ de $T$ et $S$, on consid\`ere l'op\'erateur markovien
$V$ associ\'e, de $L^2(X)$ vers $L^2(Y)$ et le facteur $\A$ de $Y$
engendr\'e par l'image de $V$, qu'on appellera le facteur engendr\'e
par $V$ ou le facteur de $Y$ engendr\'e par $\la$.

C'est aussi le facteur qu'on obtient dans la construction d'un facteur
de $Y$ isomorphe \`a un facteur d'un autocouplage infini de $X$
(voir par ex. l'article \cite{GAG}). En effet, soit $\la_{\infty}$ le produit relativement
ind\'ependant d'une infinit\'e de copies de $\lambda$ au-dessus de
leur facteur commun $Y$, sur $X^{\infty}\times Y$ et soit
$\til\la_{\infty}$ sa projection sur $X^{\infty}$. Pour des fonctions
mesurables born\'ees $f_1$, ..., $f_n$ sur $X$ et $g$ sur $Y$,
\begin{equation}\label{un}
\int \otimes_{i=1}^n f_i\otimes g\,d\la_{\infty}=
\int \bigl(\prod_{i=1}^n Vf_i\bigr)\,g\,d\nu.
\end{equation}
Comme une fonction mesurable sur $X^{\infty}\times Y$ sym\'etrique par
rapport aux variables dans $X$ est \'egale
$\la_{\infty}$-presque-partout \`a une fonction de la variable dans
$Y$ seule (``loi de $0$-$1$ relative"), le facteur sym\'etrique de
$(X^{\infty}, \til\la_{\infty})$ s'identifie \`a un facteur de $Y$, et
il est clair d'apr\`es (\ref{un}), que $\A$ est le facteur engendr\'e par
l'image de $V$ dans $L^2(Y)$.

\vskip 1ex
Dans la d\'emonstration qui suit, on applique cela \`a un 
autocouplage $\la$ de $T$. Pour \'eviter les confusions, les diff\'erents 
facteurs $X$ dans $X^{\infty}$ sont
not\'es $X_1$, ..., $X_n$,... Les projections de $\la_{\infty}$
sur $X_1\times \ldots\times X_n$ et sur $X_1\times \ldots\times
X_n\times X$ sont not\'ees $\til\la_n$ et $\la_n$ respectivement (en 
particulier $\la_1=\la$).

\begin{theoreme} Le facteur engendr\'e par un autocouplage limite de 
graphes de puissances est disjoint de tout automorphisme m\'elangeant.
En particulier, tout automorphisme non m\'elangeant a un facteur disjoint 
de tout automorphisme m\'elangeant.
\end{theoreme}
{\em D\'emonstration.\ } On consid\`ere un autocouplage $\la$ de $T$
limite faible d'une suite de graphes de puissances $T^{k_j}$ (o\`u
$k_j$ tend vers l'infini), l'op\'erateur markovien $V$ associ\'e et le facteur
$\A$ engendr\'e, avec les notations ci-dessus. Soit $(Y,\C,\nu, S)$ un
syst\`eme m\'elangeant et $\eta$ un couplage de $T$ restreint \`a $\A$
avec $S$. On
veut montrer que $\eta$ est trivial. On peut sans perte de généralité supposer $\A=\B$ (mod $\mu$).

Puisque $\la_{\infty}$ identifie $\A$ \`a un facteur de
$(X^{\infty},\til\la_{\infty})$, il revient au m\^eme de montrer la
m\^eme chose pour le couplage produit $\tilde\rho=\eta\circ\la_{\infty}$,
projection sur $X^{\infty}\times Y$ du produit relativement
ind\'ependant $\rho$ de $\la_{\infty}$ et de $\eta$ au-dessus du facteur
commun $X$.

Pour des fonctions born\'ees  $f_1$, ..., $f_n$ sur $X$ et $g$ sur $Y$, 
on a
\begin{equation}\label{deux}
\int \otimes_{i=1}^n f_i\otimes g\,d\tilde\rho=
\int \bigl(\prod_{i=1}^n Vf_i\bigr)\,J^{\ast}g\,d\mu,
\end{equation}
o\`u $J$ est l'op\'erateur markovien associ\'e \`a $\eta$, de $L^2(X)$
vers $L^2(Y)$. On doit montrer que cette int\'egrale est \'egale \`a
$$\left(\int \otimes_{i=1}^n f_i\,d\til\la_{\infty}\right)\bigl(\int
g\,d\nu\bigr)$$
et il suffira de montrer qu'elle est nulle chaque
fois que $\int g\,d\nu=0$.

Pour $n\ge 1$, notons $\rho_n$ le produit relativement ind\'ependant de $\la_n$ 
et de $\eta$ au-dessus de $X$, et $\tilde\rho_n=\eta\circ\la_n$ la projection de $\tilde\rho_n$ sur 
$X_1\times\ldots\times X_n\times Y$.  Il s'agit donc de montrer que que $\tilde\rho_n$ est le produit 
$\til\la_n\otimes\nu$ pour tout $n\ge 1$.

Cela est vrai pour $n=1$ du fait que $JV=\lim J T^{k_j}=\lim S^{k_j} J$ et $\lim  S^{k_j}$ est le projecteur 
orthogonal sur les constantes.

Fixons $n\ge1$, supposons la propri\'et\'e vraie pour cet entier, et consid\'erons 
des fonctions born\'ees $f_1$, ..., $f_{n+1}$  sur $X$ et $g$ sur $Y$, avec 
$g$ d'int\'egrale nulle. On a aussi
\begin{equation}\label{trois}
\int \otimes_{i=1}^{n+1} f_i\otimes g\,d\tilde\rho_{n+1}=
\int Vf_{n+1}\otimes\bigl(\otimes_{i=1}^n f_i\bigr)\otimes g\,d\rho_n.
\end{equation}

D'apr\`es l'hypoth\`ese de r\'ecurrence, la mesure 
spectrale de $\bigl(\otimes_{i=1}^n f_i\bigr)\otimes g$ est la 
convolu\'ee de celles des deux termes et en particulier, comme $S$ 
est m\'elangeant et $g$ d'int\'egrale nulle, sa transform\'ee de 
Fourier tend vers $0$ \`a l'infini. De l'autre c\^ot\'e, du fait que 
$V^{\ast}$ est aussi limite de puissances de $T$, et que $\Ker 
V^{\ast}=(\Im V)^{\perp}$, 
on a que la mesure spectrale de $Vf_{n+1}$ est \'etrang\`ere \`a toute 
mesure dont la transform\'ee de Fourier tend vers $0$ \`a l'infini.

Il en r\'esulte que l'int\'egrale (\ref{trois}) est nulle, ce qui 
montre la propri\'et\'e par r\'ecurrence.

\null\hfill$\blacksquare$

{\em Variante.\ }
Il n'est pas réellement nécessaire de considérer l'autocouplage infini. Voici une variante de la démonstration ne faisant appel qu'aux produits finis. On garde les mêmes notations et l'hypothèse supplémentaire $\A=\B \mod \mu$.

Comme alors $L^2(\B,\mu)$ est linéairement engendré par les produits $\prod_{i=1}^n Vf_i$ où $f_1$, ..., $f_n$ sont des fonctions
mesurables born\'ees sur $X$, on a à montrer que pour tout $n\ge 1$, étant données $n$ fonctions mesurables bornées $f_1$, ..., $f_n$ sur $X$ et une fonction mesurable bornée $g$ sur $Y$,
\begin{equation}\label{trois}
\int \bigl(\prod_{i=1}^n Vf_i\bigr)\otimes g\,d\eta=\int \bigl(\prod_{i=1}^n Vf_i\bigr) d\mu\int g\,d\nu.
\end{equation}
Pour tout $n\ge 1$ le couplage $\rho_n$ sur $X\times X^n\times Y$ produit relativement indépendant au-dessus de $X$ de $n$ copies de $\lambda$ et de $\eta$ est donné par
\begin{equation}\label{quatre}
\int f\otimes\bigl(\otimes_{i=1}^n f_i\bigr)\otimes g\,d\rho_n=\int f\cdot \bigl(\prod_{i=1}^n Vf_i\bigr)\cdot Jg\,d\mu=\int f\cdot \bigl(\prod_{i=1}^n Vf_i\bigr)\otimes g\,d\eta.
\end{equation}

En comparant (\ref{trois}) et (\ref{quatre}) lorsque $f=1$, on voit qu'on a à montrer que pour tout $n\ge 1$ les facteurs $X^n$ et $Y$ de ce couplage sont indépendants, autrement dit que la projection $\tilde\rho_n$ de $\rho_n$ sur $X^n\times Y$  est la mesure produit, et il suffit de montrer que cette intégrale est nulle lorsque $\int g \,d\nu=0$. 

Pour $n=1$ cela résulte du fait que $S$ est mélangeante et que 
\begin{align*}
\int Vf_1\cdot Jg \,d\mu=\lim \int T^{k_j}f_1\cdot Jg \,d\mu&=\lim \int T^{k_j}f_1\otimes g \,d\eta\\
&= \lim \int f_1\otimes S^{-k_j}g \,d\eta=\lim \int J^{\ast}f_1\cdot S^{-k_j}g \,d\nu.
\end{align*}
Fixons $n\ge1$, supposons la propri\'et\'e vraie pour cet entier, et consid\'erons 
des fonctions born\'ees $f_1$, ..., $f_{n+1}$  sur $X$ et $g$ sur $Y$, avec 
$g$ d'int\'egrale nulle. En appliquant (\ref{deux}) avec $f=Vf_{n+1}$, on obtient
\begin{align*}\label{trois}
\int \bigl(\prod_{i=1}^{n+1} Vf_i\bigr)\otimes g\,d\eta&=
\int Vf_{n+1}\otimes \bigl(\otimes_{i=1}^n f_i\bigr) \otimes g\,d\rho_n\\
&=\lim \int T^{k_j}f_{n+1}\otimes \bigl(\otimes_{i=1}^n f_i\bigr) \otimes g\,d\rho_n\\
&=\lim \int f_{n+1}\otimes \bigl(\otimes_{i=1}^n T^{-k_j}f_i\bigr) \otimes S^{-k_j}g\,d\rho_n\\
&=\lim \int R f_{n+1}\cdot\left(\bigl(\otimes_{i=1}^n T^{-k_j}f_i\bigr) \otimes S^{-k_j}g\right)\, d\tilde\rho_n,
\end{align*}
 où $R$ est l'opérateur markovien associé de $L^2(X, \mu)$ vers $L^2(X^n\times Y, \tilde\rho_n)$, et il suffit de montrer que $\bigl(\otimes_{i=1}^n T^{-k_j}f_i\bigr) \otimes S^{-k_j}g$ tend faiblement vers $0$. 

Or c'est une conséquence immédiate de l'hypothèse de récurrence -- les facteurs $X^n$ et $Y$ sont indépendants pour $\pi\rho_n$ -- et du fait que $S$ est mélangeant.

Il en r\'esulte que $X^{n+1}$ et $Y$ sont indépendants pour $\tilde\rho_{n+1}$, ce qui montre la propri\'et\'e cherchée par r\'ecurrence.

\null\hfill$\blacksquare$

\end{document}